

\baselineskip=14pt
\parskip=10pt

\font\eightrm=cmr8 
\font\eighttt=cmtt8
\magnification=\magstephalf

\def\1{{\overline{1}}}
\def\2{{\overline{2}}}
\parindent=0pt
\overfullrule=0in

\def\frac#1#2{{#1 \over #2}}
\bf
\centerline
{
A Eulogy for Jack Good
}
\rm
\bigskip
\centerline{ {\it
Doron 
ZEILBERGER}\footnote{$^1$}
{\eightrm  \raggedright
Department of Mathematics, Rutgers University (New Brunswick),
Hill Center-Busch Campus, 110 Frelinghuysen Rd., Piscataway,
NJ 08854-8019, USA.
{\eighttt zeilberg  at math dot rutgers dot edu} ,
\hfill \break
{\eighttt http://www.math.rutgers.edu/\~{}zeilberg/} .
Dec. 3, 2009.
Expanded version of a talk delivered Oct. 25, 2009, 3:00-3:23pm DST,
at the Algebraic Combinatorics Special Session, the 1052nd American Mathematical Society (sectional Fall) meeting,
Pennsylvania State University, Thomas Building, Room 216.
Accompanied by Maple package {\eighttt JACK}
downloadable from  the webpage of this article \hfill \break
{\eighttt http://www.math.rutgers.edu/\~{}zeilberg/mamarim/mamarimhtml/jack.html},
where one can also find sample input and output.
(The package can also be downloaded directly from \hfill\break
{\eighttt http://www.math.rutgers.edu/\~{}zeilberg/tokhniot/JACK}) .
This article is
exclusively published in the Personal Journal of Ekhad and Zeilberger 
{\eighttt http://www.math.rutgers.edu/\~{}zeilberg/pj.html}
and arxiv.org .
Supported in part by the United States of America National Science Foundation.
}
}

{\bf Kadish}

{\bf Irving John Good}, what a good name! So high-class sounding, ``waspy'', gentrified, and gentile.
I have always envisioned Jack Good- one of my {\bf greatest} heroes
and {\it influencers}- who studied in Cambridge University, and worked at {\it Bletchley Park} along
with Alan Turing and other giants, as the prototypical English (gentile) gentleman.
Imagine my surprise, when I finally met him, on Nov. 14, 2004 (when he was almost 88-years-old)
at the nice Shabbat dinner at the house of Gail Letzter and Daniel Farkas, after my colloquium talk
at Virginia Tech, and he told me that he was born {\bf Isadore Jacob Gudak}. Now {\it this} name
is as Jewish as it gets! No wonder he had to change his name, since in those days, and perhaps
even today, you can't get ahead in life with such a name, {\it even} if you have the genius of Jack Good.

Good also told me that his father, who escaped to England from the stifling {\it shtetl}
and mandatory military service, starting out as a watchmaker's apprentice,
later becoming the owner-manager of a fashionable jewelry shop in London,
was a notable Yiddish writer who published under the pen-name of {\it Mosheh Oved}. And indeed the
name rang a bell, and I am sure that I have heard that name before. Later I got a hold
of a translation of Oved's masterpiece {\it Visions and Jewels} (Faber and Faber, London, 1952) that is a very unusual
autobiography/memoir (``deserved to be read for its strangeness'' according to the dust-jacket blurb).
Part of that book describes Jack's father's adventures, starting out with his
genealogy, and it turned out that Jack's paternal grandfather was a {\it shochet} (ritual slaughterer),
that in the Jewish pecking order of the shtetl is just one notch below a rabbi (and, as it turned out,
Jack's grandfather would have become a rabbi if not for a fire that destroyed all their wealth).

Of course, already Jack's father was not observant, and definitely Jack himself was as secular as one can get,
but I am sure that he would not object to have someone say {\it kadish} for him. Unfortunately, Jack
never married, and didn't have any children, so there was no one to do it. Since I feel like
his {\it spiritual} son, let me take this opportunity to say {\it kadish} for the elevation of his soul.

[I recite, from memory, the mourner's {\it kadish}:

{\it
yitgadal ve-yitkadash shmay raba, be-alma dee bra ci-reutei, veyamlich malchutei be-khayechon u-be-yomechon
u-bekhayei de-col beit yisrael. be-agala u-bizman kariv ve-imru amen.

yehe shme raba mevorach le-olam u-le-olmei almaya .

yitbarach ve-yishtabach, ve-yitpa-ar ve-yitromam ve-yitnashe ve-yithadar ve-yitale ve-yithalal shme dikudsha brich hoo
le-ela min col birchata ve-shirata, tushbekhata ve-nekhematam da-amiran be-alma ve-imru amen.

titkabal tslothon u-ba-oothon de-col beit yisrael kedam avuhun di bishmaya ve-imru amen.

yehe shlama raba min shmaya ve-khayim tovim aleinu ve-al col yisrael ve-imru amen.
ose shalom bi-meromav hoo ya-ashe shalom aleynu ve-al col yisrael, ve-imru amen.
yehe shlama raba min shmaya ve-khayim aleinu ve-al col yisrael ve-imru amen.
}]

Jack, as a boy, is mentioned in his father's book [I wave the book
{\it Vision and Jewels} and read from it]. On p. 158 it says:

``{\it On arriving home there ran towards me, to meet me, our fiery little boy, his little thumbs
stuck into the arm-holes of his waistcoat, like a national rabbi. Pulling out his
elastic braces, thrusting his little chest forward, like Ramsay MacDonald, he said:

\qquad\qquad `Daddy, I will yet be the Honourable Master Isadore Good.' $\dots$ ''}

Well, Jack never became ``Honorable'' (in the sense of a politician or a judge), and he stopped being Isadore, but
he did become one of the greatest probabilists and statisticians of his time, and a pioneer of the once
fringe {\it Bayesian} approach to statistics. That Bayesian methods are now so accepted is due, in so small part,
to his preachings.

But I only learned of Jack's seminal work on the foundation of probability much later.
My first encounter with Jack's work was circa 1977, when Dick Askey challenged me to
prove the so-called Andrews' $q$-Dyson conjecture, that lead me to look up proofs of the original
Dyson conjecture, and inevitably to Jack Good's {\it classic} proof that {\bf changed my mathematical life}.

{\bf The Best 1-Page Mathematical Proof of All time}

The best way to eulogize a mathematician is to recite one of his proofs.
Better still, I will now give you each a copy that you should keep in your wallet.
Whenever you are feeling down, you should look at it, and I am sure that
you would cheer up.

[I now distribute to every person in the audience a copy of Good's half-page proof].

Let me reproduce it here in its {\it entirety}. 

\vfill\eject

\centerline
{Journal of Mathematical Physics \qquad Volume 11, Number 6 \qquad June 1970}

\centerline
{\bf Short Proof of a Conjecture by Dyson}

\centerline
{I.J. Good}

\centerline
{\it Department of Statistics, Virginia Polytechnic Institute, Blacksburg, Virginia}

\centerline
{(Received 26 December 1969)}

{\eightrm Dyson made a mathematical conjecture in his work on the distribution of energy levels
in complex systems. A proof is given, which is much shorter than two that have been published before.}

Let $G({\bf a})$ denote the constant term in the expansion of
$$
F( {\bf x} ; {\bf a})= \prod_{ i \neq j} \left ( 1- \frac{x_i}{x_j} \right )^{a_j} \quad , \quad i , j=1,2, \dots, n \quad ,
$$
where $a_1,a_2, \dots, a_n$ are nonnegative integers and where $F({\bf x} ; {\bf a})$ is expanded in positive and
negative powers of $x_1, \dots , x_n$. Dyson\footnote{${}^1$}
{\eightrm  \raggedright F.J.Dyson, J. Math. Phys. {\bf 3}, 140, 157, 166 (1962)}
conjectured that $G({\bf a})=M({\bf a})$, where $M({\bf a})$ is the multinomial coefficient $(a_1+ \dots + a_n)!/(a_1! \cdots a_n!).$
This was proved by Gunson\footnote{${}^2$}
{\eightrm  \raggedright J. Gunson, J. Math. Phys. {\bf 3}, 752 (1962).} and by
Wilson\footnote{${}^3$}{\eightrm  \raggedright K.G. Wilson, J. Math. Phys. {\bf 3}, 1040 (1962)}. A much shorter proof is given here.

By applying Lagrange's interpolation formula (see, for example, Kopal\footnote{${}^4$}
{\eightrm 
{\it Numerical Analysis}
(Chapman and Hall, London, 1955), p. 21}) to the function of $x$ that is identically equal to $1$ and then putting
$x=0$, we see that
$$
\sum_{j} \prod_{i} \left ( 1- \frac{x_j}{x_i} \right )^{-1} =1 \quad , \quad i \neq j \, .
$$
By multiplying $F({\bf x}; {\bf a})$ by this function we see that, if $a_j \neq 0, j=1, \dots, n$, then
$$
F({\bf x}; {\bf a})=\sum_{j} F({\bf x}; a_1, a_2, \dots, a_{j-1}, a_j -1, a_{j+1}, \dots , a_n ) \quad ,
$$
so that
$$
G({\bf a})=\sum_j G(a_1, \dots, a_{j-1}, a_j-1, a_{j+1}, \dots , a_n) \, .
\eqno(1)
$$
If $a_j=0$, then $x_j$ occurs only to negative powers in $F({\bf x};  {\bf a})$ so that
$G({\bf a})$ is then equal to the constant term in
$$
F(x_1, \dots, x_{j-1},x_{j+1}, \dots, x_n ; a_1, \dots, a_{j-1}, a_{j+1}, \dots , a_n) \, ,
$$
that is
$$
G({\bf a})=G(a_1, \dots, a_{j-1}, a_{j+1}, \dots, a_n) \, , \,\, if \quad a_j=0 \, .
\eqno(2)
$$
Also, of course
$$
G({\bf 0})=1 \, .
\eqno(3)
$$
Equations (1)-(3) clearly uniquely define $G({\bf a})$ recursively. Moreover, they are satisfied
by putting $G({\bf a})=M({\bf a})$. Therefore $G({\bf a})=M({\bf a})$, as conjectured by Dyson.

\vfill\eject

{\bf Jack Good: A very Good Man}

The brilliant ideas in Good's proof immediately lead to my article
{\it
The algebra of linear partial difference operators and its applications}, published in SIAM J. Math. Anal. 11, 919-934 (1980). 
I was very proud of that paper, and sent it out to several people, including Jack Good.
In my  youthful naivet\'e, I asked four experts, after reading the article and being
duly impressed, to pick up the phone, and call up the chair of the Illinois Math Department,
and tell him to give me a job. I needed that job desperately, since my then-girlfriend (and now wife, Jane)
was a graduate student at Urbana-Champaign.
Three out of the four politely replied that they are very busy etc., and
more to the point, they don't know me!
(one of them (John Riordan)  wrote: ``{\it you can't sell a pig in a poke}'').
But only one of the four solicited experts , {\bf Jack Good}, actually  did call!
Now that I am much older, and possibly a bit wiser, I realize how preposterous my request was,
and how kind it was on Jack's part to pick up the phone. As it turned out, 
the chairperson, Paul Bateman, refused to give me that
job that year (1978-1979), but he finally came around-in large part thanks to Jack-the following year,
and I was able to reunite with Jane, who by then became my wife.

{\bf Jack Good: A Visionary prophet of (strong!) AI}

Jack's proof of Dyson's conjecture, that I have just distributed to you, is a masterpiece of
{\bf terseness}. I experimented with it, trying to see if I can delete any word without
ruining it. I failed. Jack's proof is not only a mathematical masterpiece, but
a {\it literary} one!

Yet there are times to be terse and there are times to be verbose. Jack was one of the most prolific scientists of 
all time, and wrote many wonderful, leisurely, essays about the philosophical foundations of probability and
the probabilistic foundations of philosophy. He was also one of the earliest proponents of artificial intelligence,
and, if you believe the wikipedia article, was a consultant to Stanley Kubrick when he made {\it 2001: A Space Odyssey}
that featured HAL.

Jack's most favorite articles were collected in the book 
{\it Good Thinking: The Foundation of Probability and Its Applications}, University of Minnesota Press, 1983.
In one of these papers (in pp. 106-116) ``{\it Dynamical Probability, Computer Chess, and the Measurement of Knowledge}''
(that originally appeared in `` Machine Intelligence 8 (E.W. Elcock and D. Michie, eds.
(Wiley 1977) 139-150) one finds the following lovely quotation(p. 106), that I whole-heartedly agree with:

\qquad \qquad `` To parody Wittgenstein, what can be said at all can be said clearly {\it and it can be programmed}.''

Two pages later, one can find an even better quote:

\qquad  \qquad ``Believing, as I did (and still do), that a machine will ultimately be able to simulate all intellectual
activities of any man ...''

{\bf But} when I look again at Jack's one-page {\it proof from the book} of Dyson's conjecture, I am not
so sure. Will a computer ever be able to come up with such a gorgeous proof? But then again,
{\bf maybe it will!}

{\bf Computerized Deconstruction of Jack Good's Lovely Human Proof}

Let's try and see how a computer (once it is suitably programmed with general purpose algorithms)
would tackle Dyson's conjecture.

The {\it crux} of Good's proof is the fact that
$F( {\bf x} ; {\bf a})$ satisfies the {\it partial linear recurrence equation} with {\it constant coefficients} 
$$
F({\bf x}; {\bf a})=\sum_{j} F({\bf x}; a_1, a_2, \dots, a_{j-1}, a_j -1, a_{j+1}, \dots , a_n ) \quad , 
$$
where the coefficients neither depend on ${\bf a}$ nor on ${\bf x}$. Now the fact that there is such
a {\it simple} recurrence is indeed a miracle, but the fact that, for any specific dimension, there
is {\it some} such linear recurrence (with constant coefficients) is guaranteed a priori, and 
a computer can find that recurrence (for small $n$) rather fast. This was the main observation
of my above-mentioned 1980 paper.

Indeed, let $R_1, \dots, R_n$ are homogeneous Laurent polynomials of degree $0$ in $n$ variables,
(or equivalently, arbitrary polynomials of $n-1$ variables) and consider the Laurent polynomial
$$
F({\bf x}; {\bf a})=\prod_{i=1}^{n} R_i({\bf x})^{a_i} \quad .
$$
Introducing the shift-operators $A_i$ in the discrete variable $a_i$ ($i=1, \dots , n$), defined by
$A_i f(a_1, \dots, a_n)= f(a_1, \dots , a_{i-1},a_i+1, a_{i+1}, \dots , a_n)$, we get that
$F$ is {\it annihilated} by the $n$ operators
$$
A_i-R_i \quad , \quad (i=1 \dots n) \quad .
$$
By using the Buchberger algorithm (Gr\"obner bases) or otherwise, the computer can eliminate all the
$x$'s and get a  pure recurrence operator
$$
P(A_1, \dots, A_n) \quad,
$$
annihilating $F({\bf x}; {\bf a})$. Since such an operator is free of the $x$'s it also annihilates
each and every coefficient, in particular the {\it constant term} (the coefficient of $x_1^0 \cdots x_n^0$).

Now, if you take random $R_i$'s, $P(A_1, \dots , A_n)$ will be usually very complicated, and the
complexity gets higher with higher dimensions. The miracle of the Dyson product was that
it turned out, that for {\it every} $n$:
$$
P(A_1, \dots , A_n)=1-\sum_{i=1}^{n} A_i^{-1} \quad .
$$
A computer can discover it (and prove it!), {\it routinely} for each specific $n$, say $1 \leq n \leq 8$,
but, at present one still needs the ``human'' identity
$$
\sum_{j} \prod_{i} \left ( 1- \frac{x_j}{x_i} \right )^{-1} =1 \quad , \quad  i \neq j \quad ,
$$
that is purely routine for any specific {\it numeric} $n$, but seems to  need a human being to
prove it for {\it symbolic} $n$ (i.e. ``all'' $n >0$). Jack Good invoked the Lagrange Interpolation Formula
(whose human proof is one-line: two polynomials of degree $< n$ that coincide in $n$ different values
must be the same). But even if you never have heard of Lagrange or Interpolation,
you can still easily prove
this identity by induction on $n$, and I am sure that a computer can be taught how
to find such a proof (in the style of the Zeilberger algorithm, but in the context of symmetric functions).

 The Maple package {\tt JACK} has programs to automatically find such recurrences for any given
set of Laurent polynomials $R_1,R_2, \dots$. 

It can be gotten directly from

{\eighttt http://www.math.rutgers.edu/\~{}zeilberg/tokhniot/JACK} \quad ,

or via a link from the webpage of this article

{\eighttt http://www.math.rutgers.edu/\~{}zeilberg/mamarim/mamarimhtml/jack.html} \quad ,

where one can also find some sample input and output. In particular, automatic 
Good-style proofs of Dyson's conjecture for $n \leq 8$, from which one can
clearly see the pattern of the recurrence.

{\bf But Jack Good {\it May} have been wrong on one point ...}

Let's go back to the above-quoted Jack's prophesy, with my {\bf added} {\it emphasis}:
  	
\qquad  \qquad ``Believing, as I did (and still do), that a machine will ultimately be able to simulate {\bf all} {\it intellectual}
activities of {\bf any} man ...''

I whole-heartedly agree with Good if you replace ``{\it intellectual}'' by ``{\it mathematical}''
but I am not so sure about ``{\it intellectual}''. Will a computer {\it ever} be able to write
so beautifully and so eloquently? Isadore Jacob Gudak was much more than a mere 
mathematician and statistician, he was a {\it true} intellectual, a visionary, and
a {\it poet}, and I estimate that the (current, Bayesian!) probability that
{\it all} {\bf his} activities would be one day simulated by machine-kind is rather low.
Of course, Bayesian probabilities are always subject to change with more evidence,
so let's wait and see (and hope!).

\end